\documentclass[12pt]{article}
\usepackage[english]{babel}
\usepackage{amssymb,amsmath,graphics}
\parskip\medskipamount
\newtheorem{theorem}{Theorem}[section]
\newtheorem{Theorem}{Main Result}

\newtheorem{lemma}[theorem]{Lemma}
\newtheorem{cor}[theorem]{Corollary}
\def\<{\langle}
\def\>{\rangle}

\newcommand{\cC}{\mathcal{C}}

\newcommand{\PG}{\mathsf{PG}}

\newcommand{\cS}{\mathcal{S}}
\newcommand{\cA}{\mathcal{A}}
\newcommand{\cP}{\mathcal{P}}

\def\qed{{\hfill\hphantom{.}\nobreak\hfill$\Box$}}

\newcommand{\cL}{\mathcal{L}}
\newcommand{\cO}{\mathcal{O}}
\newcommand{\cB}{\mathcal{B}}
\newcommand{\cV}{\mathcal{V}}

\newcommand{\K}{\mathbb{K}}

\newcommand{\inc}{\mbox{\tt I}}
\newcommand{\ssW}{\mathsf{W}}
\begin{document}

\author{Koen Struyve\thanks{The first author is supported by  the Fund for Scientific Research --
Flanders (FWO - Vlaanderen)} \and Hendrik Van Maldeghem\thanks{The second author is partly supported by a
Research Grant of the Fund for Scientific Research -- Flanders (FWO - Vlaanderen)} }
\title{\bf Moufang Quadrangles of Mixed type}

\maketitle

\begin{abstract}
In this paper, we present some geometric characterizations of the Moufang quadrangles of mixed type, i.e., the
Moufang quadrangles all the points and lines of which are regular. Roughly, we classify generalized quadrangles
with enough (to be made precise) regular points and lines with the property that the dual net associated to the
regular points satisfy the Axiom of Veblen-Young, or a very weak version of the Axiom of Desargues. As an
application we obtain a geometric characterization and axiomatization of the generalized inversive planes
arising from the Suzuki-Tits ovoids related to a polarity in a mixed quadrangle. In the perfect case this gives
rise to a characterization with one axiom less than in an already known result by the second author.
\end{abstract}

\section{Introduction}

In 1974, Jacques Tits \cite{Tit:74} introduced what he called \emph{groups of mixed type}, as a certain
generalization of algebraic groups. This was motivated by the fact that certain spherical buildings arise from
such groups and Tits classified all spherical buildings of rank at least three in \cite{Tit:74}.

Roughly, the groups of mixed type of rank 2 arise when the weight of the edge of the rank 2 Coxeter diagram is
equal to the characteristic of the underlying field. Indeed, in the commutation relation of the root groups,
the weight $w$ of the edge turns up as a coefficient, and as a power. If the corresponding term does not vanish
(i.e., if in the underlying field $w$ is not equal to $0$), then we are in the generic case where we are able
to distinguish long and short roots (by the commutation relations, but also by the geometry of the
corresponding building). However, if $w=0$, i.e., if the characteristic of the underlying field is equal to $w$
(if the diagram is included in a rank 3 diagram, then only the cases $w\in\{1,2,3\}$ occur), then the
commutation relations become much more symmetric, allowing for diagram automorphisms. If the field is perfect,
not much extra happens since the symmetry is then \emph{up to the field Frobenius automorphism  $x\mapsto
x^w$}, and we only obtain an extra group automorphism (diagram automorphism). However, if the field is not
perfect, then this ``duality'' is not surjective anymore, and we obtain the peculiar situation in which the
rank 2 geometry looks symmetric, but isn't. Technically, the duality maps the geometry \emph{into} itself, but
not \emph{onto}. In other words, the geometry (building) is isomorphic to the dual of a subgeometry. On the
algebraic level, we obtain an infinite descending chain of algebraic structures, each one containing in the
next one, and the first one parameterizing the chambers in a certain panel. Since we have two different type of
panels, we have two such chains (which are mapped onto each other by the duality). The peculiar thing is now
that ``interlacing'' subchains define subgeometries and the corresponding automorphism groups are the groups of
mixed type. If the original chains consist of fields, then the interlacing chains may consist of fields, too,
but also of vector spaces. The latter only happens for $w=2$.

In the present paper, we study the case $w=2$ in a geometric way.  This is the case where the Coxeter diagram
has a weight two edge, hence a double bond.  Geometrically, this is the case of the (Moufang) generalized
quadrangles. In the (algebraically) split  case, we have a symplectic quadrangle over some field $\K$. If $\K$
has characteristic $2$ and is perfect, then this generalized quadrangle, denoted by $\ssW(\K)$, is self-dual.
If $\K$ has characteristic $2$ and is not perfect, then we are in the mixed case. There are two types of panels
here, and hence two different parameterizations. Any point row is parameterized by $\K\cup\{\infty\}$, while
any line pencil is parameterized by $\K^2\cup\{\infty\}$ (here, $\K^2$ is the field of squares of $\K$). We
obtain two chains $\K\supseteq\K^2\supseteq \K^4\supseteq\cdots$ and
$\K^2\supseteq\K^4\supseteq\K^8\supseteq\cdots$. An interlacing chain may look like
${\K'}\supseteq{\K'}^2\supseteq{\K'}^4\supseteq\cdots$, with $\K'$ a field satisfying
$\K^2\subseteq\K'\subseteq\K$.   But we may also substitute in the first chain $\K$ by a vector space $L$ over
$\K'$ contained in $\K$, and in the second chain $\K'$ by a vector space $L'$ over $\K^2$ contained in $\K'$.
This is the most general case that can occur. We denote the corresponding (Moufang) quadrangle by
$\ssW(\K,\K';L,L')$.

The quadrangle $\ssW(\K,\K';L,L')$ has an interesting geometric property.  Indeed, all its points and lines are
\emph{regular} (for precise definitions, see below).  Moreover, the dual nets associated to the regular
elements also satisfy some regularity properties. In a very weak form one can say that these dual nets satisfy
a certain Little Desargues Axiom. We will show that this axiom, together with the regularity of points and
lines, characterizes all quadrangles of mixed type. In order to answer the question of the geometric difference
between the cases where both~/~exactly one~/~none of $L$ and $L'$ are fields, we consider the Veblen \& Young
Axiom in these dual nets. We will show that, if a generalized quadrangle has enough regular points and lines,
and if the dual nets related to the regular points satisfy the Axiom of Veblen \& Young, then the quadrangle is
of mixed type and $L'$ is a field.

But also the above chains give rise to a beautiful geometric characterization of the mixed quadrangles. The
relevant property is here that a mixed quadrangle $\Gamma$ contains a subquadrangle $\Gamma'$ isomorphic to
$\Gamma$ with the property that every triad of points (a triad being a set of three pairwise noncollinear
points) in $\Gamma'$ has a center in $\Gamma$. It turns out that this property, together with the regularity of
all points, characterizes the class of mixed and symplectic quadrangles.

Another feature of the mixed quadrangles is that certain of them admit \emph{polarities}, i.e., dualities of
order $2$. In this case, the centralizer of that polarity in the little projective group of the quadrangle is a
(generalized) Suzuki group. The set of elements fixed under a polarity can be structured to a geometry which is
called a \emph{generalized inversive plane} in \cite{hvm2}. The main result of \cite{hvm2} says that the
automorphism group of these generalized inversive planes are essentially the (generalized ) Suzuki groups. In
the present paper, we use the above characterizations of the mixed quadrangles to axiomatize the generalized
inversive planes corresponding to the generalized Suzuki groups. In the perfect case, this has already been
done by the second author in \cite{hvm0}. So we relax the axioms of \cite{hvm0} to deal with the more general
cases of imperfect fields (this uses the Veblen \& Young Axiom) and vector spaces (this uses the Little
Desargues Axiom). As a corollary, these new results let us reduce the characterization in \cite{hvm0} by one axiom.

We end this introduction by mentioning that all our results hold in both the infinite and finite case. But in
the finite case there are no proper mixed quadrangles since a finite field is always perfect. All the results
of the present paper that are also valid in this improper mixed case are actually well known for finite
quadrangles. But some of our proofs give rise to alternative arguments. As an example we mention that
Theorem~\ref{theorem2} immediately implies that, if  a finite generalized quadrangle of order $q$ has an ovoid
of regular points, then all corresponding projective planes are classical.

\section{Preliminaries and Main Results}\label{MR}

\subsection{Abstract Generalized Quadrangles}
A \emph{generalized quadrangle} $\Gamma=(\cP,\cL,\inc)$ consists of a point set $\cP$, a line set $\cL$
(disjoint from $\cP$) and a symmetric incidence relation $\inc$ between $\cP$ and $\cL$ such that

\begin{itemize}
\item[(PL)] no pair of points is incident with a pair of lines and every element is incident with at least two
elements;

\item[(GQ)] for every point $x$ and every line $L$ not incident with $x$, there exists a unique point $y$ and a
unique line $M$ such that $x\inc M\inc y\inc L$.
\end{itemize}

If every element is incident with at least three elements, then the generalized quadrangle is called
\emph{thick}. As a consequence, there are a constant number of points on each line, say $s+1$ (possibly
infinite), and a constant number of lines through a point, say $t+1$. If this is the case, we say that $\Gamma$
has order $(s,t)$. Note that we adopt common linguistic expressions such as \emph{points lie on a line, lines
go through points} to describe incidence.  We will also use the notions of \emph{collinear points} and
\emph{concurrent lines} for points that are incident with a common line and lines that go through a common
point, respectively. Noncollinear points and nonconcurrent lines will also be called \emph{opposite}. If
$x\in\cP\cup\cL$ is collinear or concurrent with $y\in\cP\cup\cL$, then we write $x\sim y$. An incident
point-line pair is called a \emph{flag}.

In \cite{hvm1}, a thick generalized quadrangle is just called a generalized quadrangle, while a not necessarily
thick one is called there a weak generalized quadrangle. In the present paper, we shall also adopt this
terminology to avoid the repetition of the word `thick' in all our statements. Henceforth, with generalized
quadrangle we mean a geometry satisfying (PL) and (GQ) and such that every element is incident with at least
three elements.

Note that, if $\Gamma=(\cP,\cL,\inc)$ is a generalized quadrangle, then also $(\cL,\cP,\inc)$ is a generalized
quadrangle. We will denote the latter by $\Gamma^D$ and call it the \emph{dual} of $\Gamma$.  The \emph{duality
principle} states that in every definition and statement, one may interchange the words `point' and `line' to
obtain a new definition or statement.

Before presenting the relevant examples, we give some further definitions for abstract generalized quadrangles.

Let $\Gamma=(\cP,\cL,\inc)$ be a generalized quadrangle and let $x$ be an arbitrary point. The set of points of
$\Gamma$ collinear with $x$ will be denoted by $x^\perp$. For a set $X\subseteq\cP$, we denote by $X^\perp$ the
set of points collinear to all points of $X$, and we abbreviate $(X^\perp)^\perp$  by $X^{\perp\perp}$. If $y$
is a point opposite $x$, then $\{x,y\}^\perp$ is called the \emph{perp} of the pair $x,y$. The \emph{span} of
the pair $x,y$ is the set $\{x,y\}^{\perp\perp}$. If every span containing $x$ is also a perp (of a different
pair of points, needless to say), then the point $x$ is called \emph{regular}. Dually one defines \emph{regular
lines}. If $x$ is a regular point, then the geometry
$\Gamma^*_x=(x^\perp\setminus\{x\},\{\{x,y\}^{\perp}:y\not\sim x\},\in\mbox{ or }\ni)$ is a \emph{dual net}
(associated to $x$), i.e., it has the property that for every point $z\in x^\perp\setminus\{x\}$ and every
\emph{block} $B=\{x,y\}^\perp$, with $y$ opposite $x$, there is a unique point $z'\in B$ not collinear with $z$
(collinearity in $\Gamma_x^*$). If $\Gamma_x^*$ is a dual affine plane, then we call $x$ a \emph{projective
point}. The motivation for this terminology is that the geometry
$\Gamma_x=(x^\perp,\{\{x,y\}^\perp:y\in\cP\},\in\mbox{ or }\ni)$ is then a projective plane, called the
\emph{perp-plane} in $x$. Projective points have nice properties. For instance,  one can easily check that $x$
is a projective point if and only if the geometry $(\cP\setminus x^\perp,
\cL\cup\{\{x,y\}^{\perp\perp}:y\not\sim x\},\inc\mbox{ or }\in\mbox{ or }\ni)$ is a generalized quadrangle if
and only if every pair of distinct perps contained in $x^\perp$ meet in a unique point.

Projective points can also be approached with triads. A \emph{triad} is a triplet of pairwise opposite points.
A center of a triad $\{x,y,z\}$ is an element of $\{x,y,z\}^\perp$. A triad is called \emph{(uni)centric} if it
has a (unique) center. Now, a regular point $x$ is projective if and only if every triad containing $x$ is
centric.

In the present paper we shall mainly work with dual nets satisfying one of the two additional assumptions. We
introduce these now.

Let $\Delta=(\cP,\cL,\inc)$ be a dual net. Noncollinear points shall be called \emph{parallel}. It is easy to
see that parallelism is an equivalence relation in $\cP$. Call the dual parallel classes of points
\emph{vertical lines} and introduce a new point $\infty$ incident with all vertical lines. This way we created
a \emph{linear space} $\overline{\Gamma}=(\overline{\cP},\overline{\cL},\overline{\inc})$ (a linear space is a
point-line geometry in which every pair of distinct points is incident with a unique common line).  Let $\cV$
be the set of all vertical lines. A \emph{triangle} is a set of three distinct elements of $\cL$, pairwise
intersecting in three different  points, which are also viewed as belonging to the triangle. Two triangles are
said to be \emph{in perspective from a point $x$} if there are three different lines through $x$ of
$\overline{\Gamma}$ each incident with a unique point of each triangle.    Consider the following two
conditions:

\begin{itemize}
\item[(LD)] For every pair of triangles, which are in perspective from the point $\infty$, and for which two
pairs of corresponding sides meet on a vertical line $V$, the third pair of corresponding sides also meets on
$V$.

\item[(VY)] If a line $L$ meets two sides of a proper triangle in two distinct points, then $L$ intersects the
third side, too.
\end{itemize}

If we want to fix and include the line $V$ of {LD} in our assumptions, we more specifically say that the dual
net satisfies (LD) \emph{with respect to the vertical line $V$}.

The letters (LD) and (VY) stand for \emph{Little Desargues} and \emph{Veblen \& Young}, respectively.

Finally we introduce some notions concerning symmetry in generalized quadrangles. In general, a
\emph{collineation} of a generalized quadrangle is a permutation of the points and of the lines preserving the
incidence relation. A point $x$ of a generalized quadrangle is called a \emph{center of symmetry} if it is
regular and if the group of collineations fixing $x^\perp$ pointwise acts transitively on the set
$\{x,y\}^{\perp\perp}\setminus\{x\}$, for some, and hence for every, point $y$ opposite $x$. The dual notion is
called an \emph{axis of symmetry}.

A \emph{duality} of a generalized quadrangle is a bijection of the point set onto the line set, together with a
bijection of the line set onto the point set, preserving the incidence relation. A generalized quadrangle is
self-dual if and only if it admits some duality. A \emph{polarity} of a generalized quadrangle is a duality of
order $2$.

Let $\rho$ be a polarity of the generalized quadrangle $\Gamma$. A point (line) $x$ of $\Gamma$ is called
\emph{absolute (with respect to $\rho$)} if $x\inc x^\rho$. A flag is \emph{absolute} if it is fixed under
$\rho$. An \emph{ovoid} of $\Gamma$ is a set of points with the property that every line is incident with
exactly one element of the ovoid. It is well known that the set of absolute points with respect to a polarity
is an ovoid of the generalized quadrangle.

\subsection{The Symplectic Quadrangles}

The prototype class of examples of generalized quadrangles is the class of \emph{symplectic quadrangles}, which
are defined as follows. Let $\rho$ be a symplectic polarity in a $3$-dimensional projective space $\PG(3,\K)$
over a field $\K$. If $\cP$ is the point set of $\PG(3,\K)$, if $\cL$ is the set of lines of $\PG(3,\K)$ fixed
under $\rho$, and if $\inc$ denotes the incidence relation in $\PG(3,\K)$, then $\ssW(\K)=(\cP,\cL,\inc)$ is a
generalized quadrangle, called the \emph{symplectic quadrangle (over $\K$)}. All the points of $\ssW(\K)$ are
regular, even projective. Conversely, Schroth \cite{Sch:92} proved that any generalized quadrangle all points
of which are projective is isomorphic to a symplectic quadrangle. In fact, Theorem 6.2.1 of \cite{hvm1} asserts
that, if all points of a generalized quadrangle $\Gamma$ are regular and at least one point is projective, then
all points are projective and $\Gamma$ is a symplectic quadrangle. The first step in the proof is to show that,
if a point $x$ of $\Gamma$ is projective, then every opposite (regular) point is also projective. We record
this step as a separate lemma for later reference.

\begin{lemma}[\cite{hvm1}]\label{projective}
Let $x,y$ be two opposite points of a generalized quadrangle $\Gamma$. If $x$ is projective and $y$ is regular,
then $y$ is projective, too.
\end{lemma}

The symplectic quadrangle has a lot of symmetry. All points of $\ssW(\K)$ are centers of symmetry. Dually, all
lines of $\ssW(\K)$ are axes of symmetry if and only if $\K$ has characteristic $2$. Also, $\ssW(\K)$ is
self-dual if and only if $\K$ is a perfect field with characteristic~$2$.  Moreover, $\ssW(\K)$ admits a
polarity if and only if there exists a \emph{Tits} automorphism $\theta:\K\longrightarrow\K:x\mapsto x^\theta$,
i.e., $(x^\theta)^\theta=x^2$, for all $x\in\K$. In this case, the set of absolute points of a polarity is the
so-called \emph{Suzuki-Tits ovoid}. Viewed as a subset of points of $\PG(3,\K)$, it is also an ovoid in the
sense of Tits \cite{Tit:62}. With each ovoid of $\PG(3,\K)$ corresponds an \emph{inversive plane}, i.e.~an
incidence structure consisting of a set of points and a set of circles, which are certain subsets of the point
set, satisfying the following axioms.

\begin{itemize}
\item [{[MP1]}] Each 3 different points are contained in exactly one circle.

\item[{[MP2]}]   For each circle $C$ and each pair of points $x,y$ with $x \in C$ and $y \notin C$, there
exists an unique circle $C'$ which contains $y$ and touches $C$ in $x$.
\end{itemize}

(``Touching'' circles are circles that meet in a unique point.) The inversive planes arising from the
Suzuki-Tits ovoids have been characterized by a set of axioms by the second author in \cite{hvm0}.  We will
generalize this result below.

We end this subsection with a description of $\ssW(\K)$ using coordinates (see \cite{hvm1}). Let
$\ssW(\K)=(\cP,\cL,\inc)$ be the symplectic quadrangle over the field $\K$. Then we may take for $\cP$ the
following set:
$$\cP=\{(\infty)\}\cup\{(a):a\in\K\}\cup\{(k,b):k,b\in\K\}\cup\{(a,l,a'):a,l,a'\in\K\},$$
and for $\cL$ the set
$$\cL=\{[\infty]\}\cup\{[k]:k\in\K\}\cup\{[a,l]:a,l\in\K\}\cup\{[k,b,k']:k,b,k'\in\K\},$$
where $\infty$ is a symbol not contained in $\K$, and where incidence is given by
$$(a,l,a')\inc[a,l]\inc(a)\inc[\infty]\inc(\infty)\inc[k]\inc(k,b)\inc[k,b,k'],$$
for all $a,a',b,k,k',l\in\K$, and
$$(a,l,a')\inc[k,b,k'] \Longleftrightarrow \left\{\begin{array}{rcl}
                                                                              a' & = & ak + b,\\
                                                                              k' & = & a^2k+l-2aa'.
                                                                              \end{array}\right.
$$

We clearly see the asymmetry if the characteristic of $\K$ is unequal to $2$. If, on the other hand, the
characteristic of $\K$ is equal to $2$, then the two above formulas are equivalent if squaring is an
automorphism, i.e., the Frobenius is surjective, implying the field is perfect.

\subsection{Mixed Quadrangles and first Main Results}

Mixed quadrangles are subquadrangles of the symplectic quadrangle $\ssW(\K)$, for $\K$ an imperfect field with
characteristic $2$ (in the other case the only (thick) subquadrangles are symplectic quadrangles over
subfields). Neither the point set nor the line set of these subquadrangles can be given with a nice set of
equations in $\PG(3,\K)$, because the corresponding collineation groups are not algebraic groups. The quickest
and most elementary way to define the mixed quadrangles is using the coordinates introduced above.

So suppose $\K$ is imperfect with characteristic $2$, and let $\K^2$ be the subfield consisting of all squares.
Let $\K'$ be a subfield with $\K^2\subseteq\K'\subseteq\K$ and let $L,L'$ be subspaces of $\K,\K'$ viewed as
vector spaces over $\K',\K^2$, respectively, with $\K^2\subseteq L'$ and $\K'\subseteq L$.   We consider the
description of $\ssW(\K)$ with coordinates as above, and we now restrict the $a,a',b$ to $L$ and the $k,k',l$
to $L'$. Then we obtain a subquadrangle that we denote by $\ssW(\K,\K;L,L')$ and call a \emph{mixed quadrangle}
(the terminology in \cite{titsweiss} mentions \emph{indifferent quadrangle}, but we prefer to name the
geometries after the groups, like the symplectic quadrangle). In order to have unique notation, we also assume
that $L$ and $L'$ generate $\K$ and $\K'$ as a ring. Note that $\ssW(\K)=\ssW(\K,\K;\K,\K)$ and that
$\ssW(\K,\K^2;\K,\K^2)$ is the dual of $\ssW(\K)$ (and this dual is isomorphic to the generalized quadrangle
arising from a nonsingular quadratic form of maximal Witt index in a five-dimensional vector space over $\K$).

It is convenient to also call $\ssW(\K)$, with $\K$ perfect with characteristic $2$ a mixed quadrangle. In this
case, we also write $\ssW(\K)=\ssW(\K,\K;\K,\K)$.

In general, the dual of $\ssW(\K,\K';L,L')$ is isomorphic to $\ssW(\K',\K^2;L',L^2)$; hence the class of mixed
quadrangles is a self-dual one. Moreover, since all points of $\ssW(\K)$ are regular, so are all points of
every mixed quadrangle, and hence so are all lines of it. A famous conjecture says that every generalized
quadrangle all elements of which are regular is isomorphic to a mixed quadrangle (in the form of a problem,
this is Problem 8 in Appendix E of \cite{hvm1}). In the finite case, generalized quadrangles all of whose
points are regular are not classified, unless one requires an additional condition on the corresponding dual
nets. In \cite{Tha-Mal:97} the condition that these dual nets satisfy the Axiom of Veblen \& Young does the
job. In the present paper we will classify all generalized quadrangles with a lot of regular points and lines,
and for which the dual nets associated to the regular points satisfy the Axiom of Veblen \& Young. Postponing a
discussion of what ``a lot'' precisely means to Section~\ref{mixed}  (see Theorem~\ref{main1}), we here state
the weakest form.

\begin{Theorem}\label{MR:VY}
A generalized quadrangle $\Gamma$ is isomorphic to some mixed quadrangle $\ssW(\K,\K';L,\K')$ if and only if
all points and lines of $\Gamma$ are regular and the dual net associated to each regular point satisfies
Condition~{\rm (VY)}.
\end{Theorem}

In order to include all mixed quadrangles, we have to appeal to Condition (LD).

\begin{Theorem}\label{MR:LD}
A generalized quadrangle $\Gamma$ is isomorphic to some mixed quadrangle $\ssW(\K,\K';L,L')$ if and only if all
points and lines of $\Gamma$ are regular and the dual net associated to each regular point satisfies
Condition~{\rm (LD)}.
\end{Theorem}

Notice that, applying duality twice, the subquadrangle $\ssW(\K^2,{\K'}^2;L^2,{L'}^2)$ of $\ssW(\K,\K';L,L')$
is isomorphic to $\ssW(\K,\K';L,L')$ itself. This observation leads to the following common characterization of
symplectic and mixed quadrangles.

\begin{Theorem}\label{MR:triad}
A generalized quadrangle $\Gamma$ is a symplectic or mixed quadrangle if and only if all its points are regular
and $\Gamma$ is isomorphic with a subquadrangle $\Gamma'$ of $\Gamma$ such that each triad of  $\Gamma'$ is
centric in $\Gamma$.
\end{Theorem}

The case $\Gamma'=\Gamma$ in the above theorem boils down exactly to Schroth's result \cite{Sch:92} mentioned
above, which we shall use in our proof (alternatively, one could use Main Result~\ref{MR:VY}, which is also a
generalization of Schroth's result, but whose proof is independent of that result).

Let us finally mention that all points of a mixed quadrangle are centers of symmetry, and all lines are axes of
symmetry. Moreover,  it follows from \cite{tent} and Theorem 21.10 in \cite{titsweiss} that, if all lines of a
generalized quadrangle $\Gamma$ are axes of symmetry, and at least one point is regular, then $\Gamma$ is a
mixed quadrangle.

\subsection{Suzuki Quadrangles and more Main Results}

It is well know, see Theorem~7.3.2 of \cite{hvm1}, that a mixed quadrangle $\ssW(\K,\K';L,L')$ admits a
polarity if and only if $\K$ admits a Tits endomorphism $\theta:\K\longrightarrow \K$ (i.e.,
$(x^\theta)^\theta=x^2$) and we can choose $\K',L,L'$ such that $\K'=\K^\theta$ and $L^\theta=L'$. Hence every
polarity in $\ssW(\K,\K';,L,L')$ is the restriction of a polarity in $\ssW(\K,\K';\K,\K')$. So the case of
$L=\K$ is a kind of principal case. Every self-polar mixed quadrangle shall be called a \emph{Suzuki
quadrangle}.

Let $\rho$ be a polarity in a Suzuki quadrangle and let $\cO$ be the set of absolute points. We define a set of
circles as follows. A circle is the set of points of $\cO$ collinear to some point not contained in $\cO$.  If
we denote the family of circles by $\cC$, then we obtain a geometry $(\cO,\cC,\in\mbox{ or }\ni)$. If $\K$ is
perfect, then this is an inversive plane (a M\"obius plane). In general, it has the following properties.

\begin{itemize}
\item[{[MP1]}] Each 3 different points are contained in at most one circle.

\item[{[MP2]}] For each circle $C$ and for every pair of points $x,y \in \mathcal{P}$ with $x \in C$ and $y
\notin C$, there exists an unique circle $C'$ which contains $y$ and touches $C$ in $x$.

\item[{[CH1]}] There exist no 3 circles which are pairwise touching in different points.

\item[{[CH2]}] For each circle $C$ and every pair of points $x,y \notin C$, we have the following three
possibilities : no circle containing $x,y$ touches $C$, one circle does or all circles do.
\end{itemize}

There are a lot of geometries that satisfy the above axioms. For instance every inversive plane obtained from
an ovoid of projective $3$-space over a field with characteristic $2$. In order to further distinguish the
geometries corresponding to the polarities in mixed quadrangles, we use the observation that each circle $C$
has a very special point, which we denote by $\partial C$ and call the \emph{gnarl} of the circle. Indeed, if
$C$ is the set of points of $\cO$ collinear with the point $x\notin\cO$, then there is a unique absolute line
incident with $x$ and hence a unique point $\partial C$ of $C$ such that the line joining $\partial C$ with $x$
is absolute. Alternatively, $\partial C$ is the unique point of $C$ incident with $x^\theta$.

The function $\partial$ has the following properties.

\begin{itemize}
\item[{[ST1]}] For each pair of points $x,y$ there exists a unique circle $C$ which contains $x$ and such that
$\partial C = y$.

\item[{[ST2]}] For each circle $C$ and point $x \notin C$, there is at most one circle $C'$ which contains both
of $x$ and $\partial C$, and such that $\partial C' \in C$.

\item[{[TR]}] Let $C$ be an arbitrary circle, and let $x,y \in C$ ($\partial C \neq x \neq y \neq \partial C$).
Let $D$ be a circle through $\partial C \neq
\partial D$. For each circle $E$ different from $C$, containing both $x$ and $\partial C$, and intersecting $D$
in two distinct points $\partial C, z$, we consider the circle $E^*$ through $z$ and touching $C$ in $\partial
C$. We also consider the circle $E^{**}$ containing $y$, touching $E$ in $\partial C$. Then $E^*\cap E^{**}$ is
contained in a circle $D'$ through $\partial C$ independent of $E$.
\end{itemize}

As mentioned before, if $\K$ is perfect, then this is an inversive plane which allows us to impose a stronger
version of [MP1].

\begin{itemize}
\item[{[MP$1'$]}] Each 3 different points are contained in exactly one circle.
\end{itemize}

The properties mentioned so far characterize the generalized inversive planes arising from polarities in mixed
quadrangles.

\begin{Theorem}\label{MR:TR}
Let $\mathcal{P}$ be a set and let $\mathcal{C}$ be a set of distinguished set of subsets of $\mathcal{P}$ all
containing at least 3 elements. Also suppose there is a map $\partial : \mathcal{C} \rightarrow \mathcal{P}$
such that $\forall C \in \mathcal{C} : \partial C \in C$. We call the elements of $\mathcal{C}$ circles and if
two of them have only one point in common, we say they touch at that point. Then
$(\mathcal{P},\mathcal{C},\partial)$ satisfies the conditions {\rm [MP1], [MP2], [CH1], [CH2], [ST1], [ST2]}
and {\rm [TR]}, if and only if $\mathcal{P}$ can be embedded in a self-polar mixed quadrangle
$\ssW(\K,\K';L,L')$ as the set of absolute points of a polarity $\rho$, the set $\cC$ corresponds to the family
of sets of absolute points collinear with a nonabsolute point, and the map $\partial$ maps a circle onto its
gnarl, i.e., $\partial C$, with $C=x^\perp\cap\cP$, is the unique point of $\cP$ incident with $x^\rho$.
\end{Theorem}

If we want to restrict to self polar mixed quadrangles of type $\ssW(\K,\K';\K,\K')$, then we may introduce the
following alternative axiom (where we call a set of points \emph{cocircular} if they belong to a common
circle).

\begin{itemize}
\item[{[F]}] Let $x$ be an arbitrary point, and let $x_1,x_2,x_3$ be three points pairwise cocircular with $x$,
but not all cocircular with $x$. If  a point $y$ is cocircular with $x$ and $x_1$, and also with $x$ en $x_2$,
but if $y,x, x_1,x_2$ are not cocircular, then  $y, x, x_3$ are cocircular.
\end{itemize}

And we will show:

\begin{Theorem}\label{MR:F}
Let $\mathcal{P}$ be a set and let $\mathcal{C}$ be a set of distinguished set of subsets of $\mathcal{P}$ all
containing at least 3 elements. Also suppose there is a map $\partial : \mathcal{C} \rightarrow \mathcal{P}$
such that $\forall C \in \mathcal{C} : \partial C \in C$. We call the elements of $\mathcal{C}$ circles and if
two of them have only one point in common, we say they touch at that point. Then
$(\mathcal{P},\mathcal{C},\partial)$ satisfies the conditions {\rm [MP1], [MP2], [CH1], [CH2], [ST1], [ST2]}
and {\rm [F]}, if and only if $\mathcal{P}$ can be embedded in a self-polar mixed quadrangle
$\ssW(\K,\K';\K,\K')$ as the set of absolute points of a polarity $\rho$, the set $\cC$ corresponds to the
family of sets of absolute points collinear with a nonabsolute point, and the map $\partial$ maps a circle onto
its gnarl, i.e., $\partial C$, with $C=x^\perp\cap\cP$, is the unique point of $\cP$ incident with $x^\rho$.
\end{Theorem}

Restricting further to the case where $\K$ is a perfect field (and thus to the symplectic quadrangles
$\ssW(\K)$) we improve upon the characterization given in \cite{hvm0}, deleting one axiom.

\begin{Theorem}\label{MR:PR}
Let $\mathcal{P}$ be a set and let $\mathcal{C}$ be a set of distinguished set of subsets of $\mathcal{P}$ all
containing at least 3 elements. Also suppose there is a map $\partial : \mathcal{C} \rightarrow \mathcal{P}$
such that $\forall C \in \mathcal{C} : \partial C \in C$. We call the elements of $\mathcal{C}$ circles and if
two of them have only one point in common, we say they touch at that point. Then
$(\mathcal{P},\mathcal{C},\partial)$ satisfies the conditions {\rm [MP$1'$], [MP2], [CH1], [CH2], [ST1]} and
{\rm [ST2]}, if and only if $\mathcal{P}$ can be embedded in a projective space $\PG(3,\K)$, for some perfect
field $\K$ of characteristic 2 admitting a Tits automorphism $\theta$, such that $\mathcal{P}$ is the set of
absolute points of a polarity of a certain symplectic quadrangle $\ssW(\K)$ in $\PG(3,\K)$ and the set of
circles of $\mathcal{P}$ is equal to the set of plane sections of $\mathcal{P}$ in $\PG(3,\K)$.
\end{Theorem}

\section{Subquadrangles with Centric Triads}

In this section we prove Main Result~\ref{MR:triad}.

If $\Gamma$ is a symplectic quadrangle, then taking $\Gamma'=\Gamma$ proves that $\Gamma$ satisfies the given
condition. Now suppose $\Gamma$ is a mixed quadrangle, say $\Gamma=\ssW(\mathbb{K},\mathbb{K}';L,L')$, for
appropriate $\mathbb{K},\mathbb{K},L$ and $L'$. Then we have

$$\ssW(\mathbb{K},\mathbb{K}';L,L')\cong\ssW(\mathbb{K}^2,{\mathbb{K}'}^2;L^2,{L'}^2)\subseteq\ssW(\K^2,\K^2;\K^2,\K^2)
=\ssW(\K^2)\subseteq\ssW(\mathbb{K},\mathbb{K}';L,L'),$$

implying that every triad in $\Gamma':=\ssW(\mathbb{K}^2,{\mathbb{K}'}^2;L^2,{L'}^2)$ is centric in $\Gamma$,
since it is already centric in $\ssW(\K^2)$.

Let us now turn to the converse. So let $\Gamma=(\cP,\cL,\inc)$ be a generalized quadrangle all of whose points
are regular, and let $\theta$ be an isomorphism from $\Gamma$ to a subquadrangle $\Gamma'$. If we consider
$\theta$ as a monomorphism from $\Gamma$ to $\Gamma$, then we can consider the direct limit $\Delta$ of the
system $$\Gamma\stackrel{\theta}{\longrightarrow}\Gamma \stackrel{\theta}{\longrightarrow} \cdots
\stackrel{\theta}{\longrightarrow} \Gamma \stackrel{\theta}{\longrightarrow} \cdots. $$ It is easy to see that
$\Delta$ is a generalized quadrangle. Also, since $\Gamma$ is regular, $\Delta$ is also regular. Indeed, we
must show, for six arbitrary points $x_1,x_2,x_3,y_1,y_2,y_3$, that $x_1,x_2\in\{y_1,y_2,y_3\}^\perp$ and
$y_1,y_2\in\{x_1,x_2,x_3\}^\perp$ implies $x_3\sim y_3$ in $\Delta$. But by the definition of direct limit, all
these points have a representative in some term of the system, hence they can be considered to be points of
$\Gamma$, all points of which are regular, and hence in which the representatives of $x_3$ and $y_3$ are
collinear. But then $x_3$ and $y_3$ are collinear in $\Delta$. A similar argument, using the fact that al
triads in $\Gamma^\theta$ are centric in $\Gamma$, shows that every triad in $\Delta$ is centric, and hence
that $\Delta$ is a symplectic quadrangle by \cite{Sch:92}. Clearly $\Gamma$ can be viewed as a subquadrangle of
$\Delta$, and  since all thick subquadrangles of a symplectic quadrangle are either mixed or symplectic
quadrangles, the assertion follows.

This completes the proof of Main Result~\ref{MR:triad}.

In order to get rid of the symplectic quadrangles, or of the symplectic quadrangles in characteristic different
from $2$, one can add assumptions as follows.

\begin{cor} \begin{itemize} \item[$(i)$]
A generalized quadrangle $\Gamma$ is a mixed quadrangle or a symplectic quadrangle in characteristic $2$ if and
only if all its points are regular, at least one line is regular, and $\Gamma$ is isomorphic with a
subquadrangle $\Gamma'$ of $\Gamma$ such that each triad of  $\Gamma'$ is centric in $\Gamma$. \item[$(ii)$] A
generalized quadrangle $\Gamma$ is a mixed quadrangle but not a symplectic one if and only if all its points
are regular but at least one point is not projective, and $\Gamma$ is isomorphic with a subquadrangle $\Gamma'$
of $\Gamma$ such that each triad of  $\Gamma'$ is centric in $\Gamma$.
\end{itemize}\end{cor}

\section{Dual nets satisfying the Axiom of Veblen \& Young}\label{veblen}

Let $\Gamma=(\cP,\cL,\inc)$ be a dual net. As before, we call the dual parallel classes of points
\emph{vertical lines} and introduce a new point $\infty$ incident with all vertical lines. This way we created
a linear space $\overline{\Gamma}=(\overline{\cP},\overline{\cL},\overline{\inc})$. If two lines $L,M$
intersect in this linear space, we write $L\sim M$.  Let $\cV$ be the set of all vertical lines. Our aim is to
prove that the Condition~(LD) follows from Condition~(VY), if there exists at least one pair of
non-intersecting lines.

So henceforth we assume that $\Gamma$ satisfies (VY), and that there are at least two non-intersecting lines in
$\Gamma$. Clearly, the latter condition is equivalent to $\Gamma$ being not a dual affine plane.

We begin with defining a projective plane for every pair of intersecting lines $L,M$. Indeed, let $L,M$ be two
intersecting lines in $\Gamma$, and let $x$ be their intersection point. Then we consider the set of lines
intersecting both $L$ and $M$ in two distinct points, together with the set of lines incident with $x$ and
meeting some line $K$ that intersects $L$ and $M$ in two distinct points. We denote that set by $\cB^*$. The
point set $\cA$ is the set of points incident with at least one element of $\cB^*$, together with $\infty$. Now
add all vertical lines to $\cB^*$ by defining $\cB=\cB^*\cup\cV$.  If we denote the restriction of
$\overline{\inc}$ still by $\overline{\inc}$ (slightly abusing notation), then we claim that
$\Delta_{L,M}=(\cA,\cB,\overline{\inc})$ is a projective plane.

Indeed, this is in fact a routine check. Let us first show that two distinct lines $X,Y$ always meet. If at
least one of $X,Y$ belongs to $\cV$, or if both $X,Y$ are incident with $x$, then this is trivial. If none of
$X,Y$ is incident with $x$, then this follows directly from (VY), as by definition both $X$ and $Y$ meet both
of $L$ and $M$. If $X$ is incident with $x$, then it intersects some line $K$ which also intersects both of $L$
and $M$ in distinct points.   Since we may assume $K\neq Y$, we may also assume that $Y,K,L$ form a proper
triangle (as otherwise $Y,K,M$ form one). Now (VY) implies that $X$ meets $Y$.

Now we show that two distinct points $y,z\in\cA$ are joined by exactly one line in $\cB$. Indeed, we clearly
may assume that neither of $y,z$ coincides with $\infty$, and that they are not incident with the same vertical
line. Hence they are incident with a unique member $X\in\cL$. We must show that $X\in\cB^*$. By definition,
$y\overline{\inc}Y\in\cB^*$ and $z\overline{\inc}Z\in\cB^*$. Suppose that $Y\inc x$. Let $K\in\cB^*$ be such
that $K$ intersects $L,M,Y$ in three different points, and suppose that $y$ is not incident with $K$. Choose an
arbitrary point $y'$ incident with $K$ and not parallel to $y$. The line $Y'$ joining $y$ and $y'$ meets both
of $L$ and $M$ by (VY). We have shown that we may assume that $Y$, and hence neither $Z$, is not incident with
$x$. Moreover, using (VY), we can arrange that $Y,Z$ do not meet on $L$ or $M$ (if they do then we may
re-choose $Y$ not incident with the intersections of $Z$ with $L$ and $M$). Then $X$ meets two sides of both of
the triangles $Y,Z,L$ and $Y,Z,M$ in distinct points, and hence (VY) implies that $X$ meets both of $L$ and
$M$. If $X$ is not incident with $x$, then $X\in\cB^*$ by definition; if $x\inc X$, then with $K\in\{Y,Z\}$, we
see that again $X\in\cB^*$.

Clearly $\Delta_{L,M}=\Delta_{L',M'}$ for $L',M'$ distinct non-vertical lines of $\Delta_{L,M}$. Hence if two
projective planes like that share two non-vertical lines, then they coincide.

If we now remove from $\Delta_{L,M}$ the point $\infty$ and the vertical lines, then we obtain a dual affine
plane. Our assumptions and the existence and uniqueness of the above constructed projective plane now implies
that the dual of $\Gamma$ is a subplane covered net in the sense of Johnson \cite{Joh:95}. It follows from the
latter paper that we can identify $\cP$ with the points of a projective space $\mathbf{P}$ minus a subspace $W$
of codimension~$2$, and $\cL$ can be identified with the lines of $\mathbf{P}$ that do not intersect $W$. Our
hypothesis that $\Gamma$ is not a dual affine plane implies that the dimension of $\mathbf{P}$ is al least $3$,
and hence it is a Desarguesian projective space.

Now if a pair of triangles is in perspective from $\infty$, and if two pairs of corresponding sides meet, then
in $\mathbf{P}$, this means that the two triangles are also in perspective from a point (because two
corresponding pairs of sides must lie in the same plane), and so by Desargues' theorem, also the third pair of
corresponding sides meets, and this intersection point is collinear with the two others. This shows (LD).

Hence we have proved the following theorem.

\begin{theorem}\label{theorem1}
A dual net which is not a dual affine plane satisfies {\rm (VY)} only if it satisfies {\rm (LD)}.
\end{theorem}

One of our crucial tools to characterize the mixed quadrangles is Property (LD) for the nets associated to the
regular points of some generalized quadrangle $\Gamma$, which we now know to hold if (VY) is satisfied for that
these dual nets, on the condition that these nets are not dual affine planes. In dual affine planes (VY) holds
trivially, but (LD) is not necessarily true. A sufficient condition for (LD) is that the corresponding
projective plane is a Moufang plane. And that is exactly what we are going to prove in case that the
generalized quadrangle contains ``enough'' projective points.

\section{Generalized Quadrangles with a lot of Projective Points}\label{s:projective}

In this section we concentrate on generalized quadrangles with a number of projective points. In fact, we only
need one projective point and a set of regular points. More precisely, let $\Gamma$ be a generalized quadrangle
and let $\cO$ be a set of regular points of $\Gamma$. We assume the following two conditions on $\cO$.

\begin{itemize}
\item[(PP)] At least one member of $\cO$ is a projective point.

\item[(TP)] If $x,y$ are opposite points of $\Gamma$, then $|\{x,y\}^\perp\cap\cO|\neq 1$.
\end{itemize}

Our aim is to prove that, under these assumptions, all points of $\cO$ are projective and every corresponding
perp-plane is a Moufang projective plane. We will need the following characterization of Moufang projective
planes proved by the second author in \cite{Mal:02}. In a projective plane, a line $L$ is called an \emph{axis
of transitivity} if the pointwise stabilizer of $L$ acts transitively on the points not incident with $L$.
\begin{theorem}[\cite{Mal:02}]\label{cohen}
A projective plane is a Moufang plane if and only if each line $L$ is an axis of transitivity.
\end{theorem}

Henceforth $\Gamma$ is a generalized quadrangle with $\cO$ a set of regular points of $\Gamma$ satisfying
(PP) and (TP).

We start with proving that all elements of $\cO$ are projective.

\begin{lemma}\label{lemma3}
Every element of $\cO$ is a projective point of $\Gamma$.
\end{lemma}

\emph{Proof.} We know that there is at least one point $p\in\cO$ which is projective. Let $q$ be any other
element of $\cO$. If $q$ is opposite $p$, then Lemma~\ref{projective} implies that $q$ is projective. Now
suppose $q\sim p$. Let $x,y$ be opposite points collinear to $p$ such that $x$ is incident with the line $pq$,
but $x\neq q$. Then $p\in\{x,y\}^\perp$, implying by (TP) that some other element $p'\in\cO\setminus\{p\}$ also
belongs to $\{x,y\}^\perp$. Clearly, $p'$ is opposite $p$ and therefore is a projective point. But $p'$ is also
 opposite $q$ and hence Lemma~\ref{projective} implies that $q$ is projective.

The lemma is proved. \qed

We now prove a lemma that will generate collineations of the perp-planes $\Gamma_p$, for $p\in\cO$.

\begin{lemma}\label{lemma4}
Let $p,q \in \mathcal{O}$, with $p$ opposite $q$. Then the following function $\theta_{p,q}$ defines an
isomorphism between $\Gamma_p$ and $\Gamma_q^D$ :
\begin{itemize}
\item[$(i)$] A point $x$ of $\Gamma_p$ is mapped to the block $x^{\theta_{p,q}}$ of $\Gamma_q$ consisting of
all the points collinear with both $x$ and $q$.

\item[$(ii)$] A block $\alpha$ of $\Gamma_p$ is mapped to the point $\alpha^{\theta_{p,q}}$ of $\Gamma_q$
collinear with $q$ and with all points of $\alpha$.
\end{itemize}
\end{lemma}

\emph{Proof.} First we show that $\theta_{p,q}$ is well defined by proving that for each block $\alpha$ of
$\Gamma_p$, there is indeed a unique point $a\sim q$ collinear with all points of $\alpha$. Indeed, we may
assume that $\alpha\neq\{p,q\}^\perp$, as otherwise $a=q$ is easily seen to be that unique point. Since
$\Gamma_p$ is projective, there is a unique point $r\in\{p,q\}^\perp\cap\alpha$. No $a$ is necessarily the
unique point on the line $aq$ which is collinear with any point of $\alpha\setminus\{r\}$.

The definition of $\theta_{p,q}$ now easily implies that, if $x\in\alpha$, with $x\sim p$ and $\alpha$ a block
of $\Gamma_p$, then $\alpha^{\theta_{p,q}}\in x^{\theta_{p,q}}$. Also, the inverse mapping is apparently given
by $\theta_{q,p}$, hence $\theta_{p,q}$ is bijective and so defines an isomorphism from $\Gamma_p$ to the dual
of $\Gamma_q$. \qed

Note that we can write $x^{\theta_{p,q}}=\{q,x\}^\perp$ and $\alpha^{\theta_{p,q}}=\alpha^{\perp\perp}\cap
q^\perp$, with $x\sim p$ and $\alpha$ a block of $\Gamma_p$.

We now consider three different points $p_1,p_2,p_3 \in \mathcal{O}$, with $p_3$ opposite both $p_1$ and $p_2$.
By the previous lemma, we can combine $\theta_{p_1,p_3}$ and $\theta_{p_3,p_2}$ to an isomorphism $\phi:=
\theta_{p_1,p_3} \theta_{p_3,p_2}$ between $\Gamma_{p_1}$ and $\Gamma_{p_2}$. Let us calculate the image of a
point $x$ of $\Gamma_{p_1}$ under $\phi$.
\begin{equation}
x^\phi = x^{\theta_{p_1,p_3} \theta_{p_3,p_2}} = (\{x,p_3 \}^\perp )^{\theta_{p_3,p_2}} = \{x,p_3 \}^{\perp
\perp} \cap p_2^\perp
\end{equation}
If we apply this to a point $a$ in $\{p_1,p_2 \}^\perp $, then, since $a\in\{a,p_3\}^{\perp\perp}\cap
p_2^\perp$, we see that $a^\phi = a$ (note the independence of $p_3$). We also have $p_1^\phi = \{p_1,p_3
\}^{\perp \perp} \cap p_2^\perp$.

Now let $p_3'$ be another point of $\mathcal{O} \backslash \{p_1,p_2\}$ opposite both $p_1,p_2$. We obtain a
different isomorphism $\phi' := \theta_{p_1,p_3'} \theta_{p_3',p_2}$ between the two perp-planes $\Gamma_{p_1}$
and $\Gamma_{p_2}$. This allows us to construct a collineation $\tau := {\phi}^{-1} \phi'$ of $\Gamma_{p_2}$.
Using the independence mentioned in the above paragraph we see that $\{p_1,p_2 \}^\perp $ is fixed pointwise
under the action of $\tau$. Choose 2 points $x,y$ in $\Gamma_{p_2}$ different from $p_2$ and not an element of $\{p_1,p_2
\}^\perp$. We can choose $p_3\in\cO$ in such a way that $p_1^\phi = x$ (this is possible since the span
$\{p_1,x\}^{\perp\perp}$ contains at least two points of $\cO$, and we can choose $p_3$ as one of them
different from $p_1$;  then $p_1^\phi = \{p_1,p_3 \}^{\perp \perp} \cap p_2^\perp = \{p_1,x \}^{\perp \perp}
\cap p_2^\perp = x$). Analogously, we can choose $p_3'\in\cO$ in such a way that $p_1^{\phi'} = y$. Combining
this we obtain $x^\tau = x^{{\phi}^{-1} \phi'} = {p_1}^{\phi'}= y$.

Consequently, the pointwise stabilizer of $\{p_1,p_2 \}^\perp $ in the collineation group of $\Gamma_{p_2}$
acts transitively on all the other points of the plane possibly except $p_2$. But if $p_2$ was fixed by this
stabilizer, then the orbits of the other points would completely lie on lines through $p_2$, which is
impossible by the transitivity already shown. So the pointwise stabilizer of $\{p_1,p_2\}^\perp$ is transitive on all
points of the perp-plane $\Gamma_2$ except for the points of $\{p_1,p_2 \}^\perp $ itself. Hence $\{p_1,p_2\}^\perp$
is an axis of transitivity in the projective plane $\Gamma_{p_2}$.

We can even do better.

\begin{lemma}\label{lemma5}
Each block $\alpha$ of $\Gamma_{p_2}$ is an axis of transitivity.
\end{lemma}

\emph{Proof.} If $\alpha$ is a perp $\{x,y\}^\perp$, then the span $\{x,y\}^{\perp\perp}$ is a perp and
contains $p_2$, hence it contains a second point $p_4\in\cO$. This implies $\alpha=\{p_2,p_4\}^\perp$ and the
assertion follows from our previous discussion.

Now the blocks through $p_2$ are also axes of transitivity because they can be mapped to blocks not through
$p_2$ by the pointwise stabilizers of the blocks not containing $p_2$, for which the condition is already true
and hence which have rich enough point stabilizers to do so. \qed

Now Theorem~\ref{cohen} implies that $\Gamma_{p_2}$, and hence all perp-planes of points in $\cO$, are Moufang
projective planes, and in particular satisfy Condition (LD).

Hence, in this section, we have shown the following theorem.

\begin{theorem}\label{theorem2}
Let $\Gamma$ be a generalized quadrangle and let $\cO$ be a subset of regular points of $\Gamma$ satisfying
{\rm (PP)} and {\rm (TP)}. Then all points of $\cO$ are projective and all corresponding perp-planes are
Moufang projective planes and satisfy in particular {\rm (LD)}.
\end{theorem}

\section{Quadrangles with Regular Points satisfying (LD)}\label{LD}

In this section, we will prove Main Result~\ref{MR:VY} and Main Result~\ref{MR:LD}. These will follow from
Theorem~\ref{theorem1} and Theorem~\ref{theorem2} and the following lemma.

\begin{lemma}\label{lemmaLD}
Let  $\Gamma=(\cP,\cL,\inc)$ be a generalized quadrangle containing a flag $\{p,L\}$ consisting of a regular
line $L$ and a regular point $p$. Then the dual net associated to $p$ satisfies {\rm (LD)} with respect to the
vertical line defined by $L$ if and only if $L$ is an axis of symmetry for $\Gamma$.
\end{lemma}

\emph{Proof.} We begin by noticing that, if there are only three lines through each point in $\Gamma$, then
$\Gamma$ has order $2$ and is isomorphic to $\ssW(2)$, in which the assertion clearly holds. So we may assume
that there are at least four lines through each point.

First we assume that the dual net associated to $p$ satisfies {\rm (LD)} with respect to the vertical line
defined by $L$

Let $M$ be a line through $p$ different from $L$. Let  $a,a'$ be two points incident with $M$ but different
from $p$. We will gradually construct a collineation $\theta$ mapping $a$ to $a'$ fixing $L$ pointwise, and
fixing all lines meeting $L$.
\subsubsection*{Lines intersecting $L$}
For these lines $N$ we set $N^\theta = N$.
\subsubsection*{Points collinear to $p$ not on $L$}
Let $N$ be a line through $p$ different from both $L$ and $M$, and let $q$ be a point on $N$ different from
$p$, then we define the image of $p$ under $\theta$ as following. The perp $\alpha$ in $\Gamma_p^*$ through
$a$ and $q$ intersects $L$ in a point $b$. Then $q^\theta$ is the intersection point of $N$ with the perp
through $a'=a^\theta$ and $b$. This way the image of $a$ defines the image of a point $q$ collinear with $p$,
but not with $a$. We denote this as : $a \rightarrow q$. The image of a point $c$ on $M$ is defined by
$q\rightarrow c$, for some point $q\sim p$ not collinear to $c$.

To show that $\theta$ is well defined we have to show that combining $a \rightarrow b$ with $b \rightarrow c$
(we will abbreviate this as $a \rightarrow b \rightarrow c$) with $b$ not collinear with either $a$ or $c$, is
independent of the choice of $b$. So suppose $a,b,c$ and $d$ are four points in $p^\perp$ not on $L$ such that
both $b$ and $d$ are not collinear with either $a$ or $c$.
\begin{itemize}
\item[$(i)$] If $a$ is not collinear with $c$ then $a \rightarrow b \rightarrow c$ is equivalent with $a
\rightarrow c$. Indeed, this follows directly from the condition (LD) applied to the triangles $a,b,c$ and
$a^\theta,b^\theta,c^\theta$ (where $\theta$ is defined using $a \rightarrow b \rightarrow c$). Similarly,
$a\rightarrow d \rightarrow c$ is equivalent with $a\rightarrow c$ and the result follows.

\item[$(ii)$] Suppose that $a$ is collinear with $c$. If $b$ is not collinear with $d$ then $a \rightarrow b
\rightarrow c$ is equivalent with $a \rightarrow b \rightarrow d \rightarrow c$ which on its turn is equivalent
with $a \rightarrow d \rightarrow c$. If $b$ and $d$ are collinear then we can choose a point $e$ collinear
with $p$ but not with $a$ or $b$ and not on $L$ (because there are at least four lines through a point in
$\Gamma$). Then $a \rightarrow b \rightarrow c$ is equivalent with $a \rightarrow b \rightarrow e \rightarrow
c$, $a \rightarrow e \rightarrow c$ and $a \rightarrow d \rightarrow c$ by using the previous arguments.
\end{itemize}
It is important to note that $\theta$ preserves the perps in $\Gamma_p^*$.

\subsubsection*{Lines and points opposite $L$ or $p$}
Let $N$ be a line opposite $L$, and let $p\inc A\inc q\inc N$. Then we define $N^\theta$ to be the unique line
in the (line) span containing $L$ and $N$ incident with $q^\theta$. The image of a point $t$ incident with $N$
is defined as the intersection point of $N^\theta$ with the unique line $K$ through $t$ intersecting $L$ (these
indeed intersect because of the regularity of $L$). The only thing left to show is that $t^\theta$ is well
defined. If $t$ is collinear with $p$ then this is clear, so suppose $t\not\sim p$. The lines through $t$
define a perp in $\Gamma_p^*$, which will be mapped to another perp by $\theta$ while fixing the intersection
point $r$ of $K$ and $L$ of the perp. The images of all the lines through $t$ must meet $K$. Since they also
must contain a point of the perp $\{p,t^\theta\}^\perp$, we see that they are all incident with $t^\theta$.
Hence $t^\theta$ is well defined. It is now also clear that $\theta$ and its inverse preserve incidence, and
hence it is a symmetry. Since $a$ and $a'$ were basically arbitrary, it follows that $L$ is an axis of
symmetry.

Now we assume that $L$ is an axis of symmetry for $\Gamma$, and that two triangles $T,T'$ in the dual net
associated to $p$ are in perspective from $p$, with two pairs $(A,A')\in T\times T'$ and $(B,B')\in T\times T'$
of corresponding sides meeting on $L$. It follows easily that the unique symmetry $\theta$ which maps $A\cap B$
to $A'\cap B'$ maps $A$ to $A'$ and maps $B$ to $B'$. If $(C,C')$ is the third pair of corresponding sides,
then, since $\theta$ preserves all vertical lines, we see that $\theta$ maps $A\cap C$ and $B\cap C$ to $A'\cap
C'$ and $B'\cap C'$, respectively. Consequently $theta$ maps $C$ to $C'$ and so, since the points of $C$ and
$C'$ on $L$ are fixed, the assertion follows. \qed

We are now ready to prove slightly more general results than Main Results~\ref{MR:VY} and~\ref{MR:LD}.

\begin{theorem}\label{main1}
A generalized quadrangle $\Gamma=(\cP,\cL,\inc)$ is a mixed quadrangle if and only if there is a nonempty
subset $\cO\subseteq\cP$ of points and a subset $\cS\subseteq\cL$ of lines satisfying the following conditions.
\begin{itemize}
\item[$(i)$] All points of $\cO$ and all lines of $\cS$ are regular.

\item[$(ii)$] Every (line) span containing a line of $\cS$ contains at least two lines of $\cS$.

\item[$(iii)$] Every element of $\cS$ is incident with some element of $\cO$.

\item[$(iv)$] The dual net associated to each regular point $x$ of $\cO$ satisfies {\rm (LD)} with respect to a
vertical line given by some element of $\cS$ incident with $x$.
\end{itemize}

In particular, if all elements of $\Gamma$ are regular and $(iv)$ holds, then $\Gamma$ is a mixed quadrangle.
\end{theorem}

\emph{Proof.} If $\Gamma$ is a mixed quadrangle, then we take for $\cO$ the point set of $\Gamma$ and for $\cS$
the line set. Since all points and lines are regular, $(i)$ up to $(iii)$ follow. Since every line is an axis
of symmetry, $(iv)$ follows from Lemma~\ref{lemmaLD}.

Now suppose the quadrangle $\Gamma$ satisfies the given properties. Fix a line $L$ of $\cS$. By $(iii)$, there
is a regular point $p$ incident with $L$ with the property that, by $(iv)$, the associated dual net satisfies
(LD). Lemma~\ref{lemmaLD} implies that $L$ is an axis of symmetry. Likewise, every element of $\cS$ is an axis
of symmetry. Let $M$ be an arbitrary line opposite $L$. The span $\{L,M\}^{\perp\perp}$ contains some element
$K\in\cS\setminus\{L\}$, by $(ii)$. Since $L$ is an axis of symmetry, there is a collineation mapping $K$ to
$M$. Since $K$ is an axis of symmetry, so is $M$. Hence all lines opposite $L$, and likewise all lines opposite
$K$, are axes of symmetry. It is easy to see that  for each element $N$ of $\{L,K\}^\perp$ there is a line
opposite all of $L,K,N$. we conclude that all lines of $\Gamma$ are axes of symmetry. Since we have at least
one regular point, we conclude that $\Gamma$ is a mixed quadrangle. \qed

\begin{theorem}\label{main2}
A generalized quadrangle $\Gamma=(\cP,\cL,\inc)$ is isomorphic to a mixed quadrangle $\ssW(\K,\K';L,\K')$ if
and only if there is a nonempty subset $\cO\subseteq\cP$ of points and a subset $\cS\subseteq\cL$ of lines
satisfying the following conditions.
\begin{itemize}
\item[$(i)$] All points of $\cO$ and all lines of $\cS$ are regular.

\item[$(ii)$] Every span containing a point of $\cO$ contains at least two points of $\cO$.

\item[$(ii)'$] Every (line) span containing a line of $\cS$ contains at least two lines of $\cS$.

\item[$(iii)$] Every element of $\cS$ is incident with some element of $\cO$.

\item[$(iv)$] The dual net associated to each regular point of $\cO$ satisfies {\rm (VY)}.
\end{itemize}

In particular, if all elements of $\Gamma$ are regular and $(iv)$ holds, then $\Gamma$ is isomorphic to a mixed
quadrangle $\ssW(\K,\K';L,\K')$.
\end{theorem}

\emph{Proof.} First suppose that the quadrangle $\Gamma$ satisfies the given properties. If none of the points
of $\cO$ are projective, then Theorem~\ref{theorem1} implies that, together with $(iv)$,  each dual net
associated to a regular point of $\cO$ satisfies (LD). From Theorem~\ref{main1} we infer that $\Gamma$ is
isomorphic to a mixed quadrangle $\ssW(\K,\K';L,L')$.  We now show that $L'=\K'$. Assume, by way of
contradiction, that $L'\neq\K'$. Then we can choose elements $k,k'\in L'$ such that $kk'\notin L'$. One easily
calculates that, in the coordinate representation of $\ssW(\K,\K';L,L')$ the perp
$T_{a,a'}:=\{(\infty),(a,l,a')\}^\perp$ consists of the point $(a)$ together with the points $(x,ax+a')$, $x\in
L'$. Now we consider the perps $T_{0,0}=\{(0)\}\cup\{(x,0):x\in L'\}$ and $T_{0,1}=\{(0)\}\cup\{(x,1):x\in
L'\}$, which both meet the perps $T_{1,0}=\{(1)\}\cup\{(x,x):x\in L'\}$ and
$T_{(k^{-1}+1)^{-1},k'(k^{-1}+1)^{-1}}=\{((k^{-1}+1)^{-1})\}\cup\{(x,(k^{-1}+1)^{-1}x+(k^{-1}+1)^{-1}k'):x\in
L'\}$. By (VY), the latter two perps must intersect. Hence there must exist $x\in L'$ such that
$$x=(k^{-1}+1)^{-1}x+(k^{-1}+1)^{-1}k',$$ which is equivalent with $kk'=x\in L'$, a contradiction.

If at least one point of $\cO$ is projective, then by Theorem~\ref{theorem2} and Assumption~$(ii)$, all points
of $\cO$ are projective, all corresponding perp-planes are Moufang and satisfy (LD). Since they also satisfy
(VY), the result now again follows from Theorem~\ref{main1} and the computation performed in the previous
paragraph.

If $\Gamma\cong\ssW(\K,\K';L,\K')$, then an elementary calculation as above shows that the dual net associated
to an arbitrary point satisfies (VY). \qed

\label{mixed}

\section{Generalized Suzuki-Tits Inversive Planes}

In this section, we generalize the main theorem of \cite{hvm0} to all self-polar mixed quadrangles. It will
turn out that we need exactly the more general form in the previous section of our Main Results~\ref{MR:VY}
and~\ref{MR:LD} in order to prove Main Results~\ref{MR:TR} and ~\ref{MR:F}..

In this section, we let $\mathcal{P}$ be a set and $\mathcal{C}$ a distinguished set of subsets of
$\mathcal{P}$ all containing at least 3 elements. Also we have been given a map $\partial : \mathcal{C}
\rightarrow \mathcal{P}$ such that $\forall C \in \mathcal{C} : \partial C \in C$. We call the elements of
$\mathcal{C}$ \emph{circles} and if two of them have only one point in common, we say they \emph{touch at that
point}. The element $\partial C$ of a circle $C$ will be called the \emph{gnarl} of $C$. We assume that
$(\mathcal{P},\mathcal{C},\partial)$ satisfies the conditions [MP1], [MP2], [CH1], [CH2], [ST1], [ST2] and
[TR].

First, we will prove some further properties using these axioms. All these lemmas are copies or reformulations
of lemmas in \cite{hvm0}, with similar proofs, although [MP1] and [ST2] are in the present paper slightly
weaker than the corresponding axioms in \cite{hvm0}. We mention them without proof.
\begin{lemma} \label{lemma:1}
Suppose we have 3 different circles $C,D$ and $E$. If $C$ and $E$ both touch $D$ at some point $x$, then $C$
touches $E$ at $x$.
\end{lemma}

\begin{lemma} \label{lemma:2}
For every circle $C$ and every point $x$ not contained in $C$ there exists an unique circle $D$  with $\partial
D \in C$, $\partial C \neq \partial D$ and containing both of $x$ and $\partial C$.
\end{lemma}

\begin{lemma} \label{lemma:3}
If a circle $C$ touches $D$ at $\partial D$, then $\partial C= \partial D$.
\end{lemma}

We now proceed with constructing a geometry $\Gamma=(\cP^*,\cL^*,\inc)$ out of
$(\mathcal{P},\mathcal{C},\partial)$. Also this is similar to the perfect case in \cite{hvm0}, but since it is
crucial for the rest, we repeat it here.

We identify both $\cP^*$ and $\cL^*$ with the union of $\mathcal{P}$ and $\mathcal{C}$. To avoid confusing the
elements of $\cP^*$ with those of  $\cL^*$, we put a subscript $p$ or $l$ to denote to which set it belongs,
i.e., for all $ x \in \mathcal{P}$ and all $C \in \mathcal{C}$, we have $x_p,C_p \in \cP$ and $x_l,C_l \in L$. A
point $x_p$, $x \in \mathcal{P}$ is incident with $y_l$, $y \in \mathcal{P}$ if and only if $x=y$. A point
$x_p$,  $x \in \mathcal{P}$ is incident with the line $C_l$, $C \in \mathcal{C}$ if and only if $C_p$ is
incident with $x_l$ if and only if $\partial C = x$. Finally, the point $C_p$, $C \in \mathcal{C}$ is incident
with $D_l$, $D \in \mathcal{C}$ if and only if $\partial C \in D$, $\partial D \in C$ and $\partial C \neq
\partial D$. This new geometry $\Gamma$ obviously admits a polarity $\rho : \cP^*\leftrightarrow\cL^*: C_p
\mapsto C_l, x_p \mapsto x_l, C_l\mapsto C_p, x_l\mapsto x_p$. The absolute flags are of the form $\{x_p,x_l\}$
with $x \in \mathcal{P}$.

The following lemma tells us when two points are collinear in $\Gamma$.
\begin{lemma} \label{lemma:coll} For all $x,y\in\cP$ and $C,D\in\cC$, the following holds.
\begin{itemize}
\item[$(i)$] The point $x_p$ is collinear with the point $y_p$ if and only if $x=y$.

\item[$(ii)$] The point $x_p$ is collinear with the point $C_p$ if and only if $x \in C$.

\item[$(iii)$] The point $C_p$ is collinear with the point $D_p$ if and only if $C$ and $D$ touch each other.
\end{itemize}
Also, two different elements of $\cP^*$ are incident with at most one element of $\cL^*$.
\end{lemma}

\emph{Proof.}
\begin{itemize}
\item[$(i)$] Suppose $x_p \inc C_l \inc y_p$, then, by definition, $x= \partial C=y$.

\item[$(ii)$] If $x_p$ is collinear with $C_p$, then $x_p \inc x_l \inc C_p$  or there is an $E \in
\mathcal{C}$ such that $x_p \inc E_l \inc C_p$. In the first case we have $x= \partial C \in C$; in the second
case $x =
\partial D \in C$. Suppose now that $x \in C$. If $x = \partial C$, then $x_p \inc x_l \inc C_p$ and so $x_p$
is collinear with $C_p$. If $x \neq \partial C$, then there is an unique circle D with gnarl $x$ through
$\delta C$ by [ST1], so $x_p \inc D_l \inc C_p$.

\item[$(iii)$] If $C_p \inc z_l \inc D_p$, with $z \in \mathcal{P}$, then the claim follows from [ST1]. Suppose
that $C_p \inc E_l \inc D_p$, with $E \in \mathcal{C}$. Then $\partial E \in C \cap D$, and since $D \neq C$,
we have $\partial D \neq \partial C$. Clearly, also $\partial C \neq \partial E \neq \partial D$. Since
$\partial C, \partial D \in E$, the result follows from [ST2].

Conversely, suppose $C$ and $D$ touch. If they touch at $\partial C$ then by Lemma~\ref{lemma:3}, $\partial C =
\partial D$ and $C_p \inc (\partial C)_l \inc D_p$. So we can assume that they touch at a point $x$ different
from $\partial C$ and different from $\partial D$. Let $E$ be the circle containing $\partial D$ and so that
$\partial E = x$, and assume by way of contradiction that $\partial C\notin E$. By Lemma~\ref{lemma:2} there
exists a circle $F$ containing $\partial C$ and $x$, and with $\partial F \in E$. Our assumption implies $F
\neq E$. We claim that either $D=F$ or $F$ touches $D$ at $x$. Indeed, if not, then $D$ and $F$ share some
point $y\neq x$. Note that $y\notin E$ as otherwise $F$ and $D$ coincide with $E$, a contradiction. But then
both $D$ and $F$ have their gnarl on $E$, contain the gnarl of $E$ and contain a further point $y\notin E$.
Lemma \ref{lemma:2} implies that $D=F$. Our claim follows. Now by lemma \ref{lemma:1}, $F$ touches $C$ at $x$,
contradicting $\partial C \in F \cap C$. So we have that $C_p \inc E_l \inc D_p$. \qed
\end{itemize}

Our goal now is to show that $\Gamma$ is a Suzuki quadrangle. First we prove that $\Gamma$ is a generalized
quadrangle.

\begin{lemma}
There are no three different, pairwise collinear points in $\cP^*$ unless they are all incident with the same
line.
\end{lemma}

\emph{Proof.} First suppose one of the points is of the form $x_p$ with $x \in \mathcal{P}$, then the other
points must be of the form $C_p$ and $D_p$ ($C,D \in \mathcal{C}$) with $x = C \cap D$. If $x=\partial C$, then
$x=\partial D$ and all the points are incident with the line $x_l$. If $x\neq \partial C$, then $C_p\inc
E_l\inc D_p$, with $E\in\cC$ and hence $\partial E=x$. But then also $x_p\inc E_l$.

Now suppose we have three points of the form $C_p,D_p$ and $E_p$ with $C,D,E \in \mathcal{C}$. By collinearity,
the circles $C,D$ and $E$ all have to touch each other. Axiom~[CH1] implies that they touch in one common point
$x$. So $C_p,D_p$ and $E_p$ are all collinear with $x_p$. By the first part of the proof we obtain that
$C_p,D_p, x_p$ lie on one line $F_l$ and $C_p, E_p, x_p$ lie one line $G_l$ ($F,G \in \mathcal{C}$). Both $F_l$
and $G_l$ contain $C_p$ and $x_p$, so, by the last assertion of Lemma~\ref{lemma:coll}, $C_p,D_p$ and $E_p$ all
are incident with $F_l = G_l$. \qed

\begin{lemma}
A point in $\cP^*$ and a line in $\cL^*$ lie at distance at most $3$ from each other.
\end{lemma}
\emph{Proof.} We prove that for any point $X$ and any line $M$ not incident with $X$, there is a point on $M$
collinear with $X$.
\begin{itemize}
\item[Case 1.] First suppose $X=x_p$ and $M=y_l$, with $x,y \in \mathcal{P}, x\neq y$. Condition~[ST1] tells us
that there is a circle C with gnarl $x$ trough $y$. Now $C_p$ is collinear with $y_p$ (by Lemma
\ref{lemma:coll}) and incident with $x_p$ (since $\partial C =x$).

\item[Case 2.] Secondly suppose $X=x_p$ and $M=C_l$, with $x \in \mathcal{P}$, $C \in \mathcal{L}$, and
$\partial C \neq x$). If $x \in C$ then the point $D_p$ with $D$ the circle with gnarl $x$ through $\partial C$
is incident with $C_l$ and collinear with $x_p$.

If $x$ is not on $C$ then by Lemma~\ref{lemma:2} there exists a circle $D$ through $x$ sharing two distinct
points (namely, $\partial C$ and $\partial D$) with $C$. The point $D_p$ is now on $C_l$ and collinear with
$x_p$.

\item[Case 3.]  Taking duality in account, there is one case left to show, where $X=C_p$ and $M=D_l$, with $C,D
\in \mathcal{L}$ and $C_p$ not incident with $D_l$ in $\Gamma$. The first possibility is that $\partial C =
\partial D$. Then $C_p$ is collinear with $(\partial C)_p$ which is incident with $D_l$.

No suppose that $\partial C \neq \partial D \in C$. Then the point $(\partial D)_p$ is collinear with $C_p$ and
lies on $D_l$. The case where $\partial C \in D$ is the dual of the case just handled.

So we may assume that $\partial C \notin D, \partial D \notin C$. By Axiom~[MP2] and the fact that a circle
contains 3 or more points, there are at least two circles circle $C_1$ and $C_2$ with gnarl $\partial D$ and
touching $C$. By Axiom~[CH1] these 2 circles have a second point $x \neq \partial D$ in common. Due to [CH2]
all circles through $x$ and $\partial D$ touch C. So we can consider the circle $E$, guaranteed to exist by
Lemma~\ref{lemma:2}, which contains the two points $\partial D, x$, and has its gnarl on $D$. This circle $E$
touches $C$; hence $E_p$ is collinear with $C_p$ and is incident with $D_l$. \qed
\end{itemize}

Now we want to apply Theorem~\ref{main2}. Hence we have to find a suitable set of regular points and regular
lines. We will consider the set of absolute points and absolute lines of $\Gamma$ with respect to the polarity
$\rho$ mentioned above.

\begin{lemma}
The absolute points and lines of $\Gamma$ are regular.
\end{lemma}
\emph{Proof.} Because of the polarity $\rho$, we only need to prove that when three different points
$\{U,V,W\}$ are collinear with two non-collinear points $X,Y$, with $X=x_p$ for some $x\in \mathcal{P}$, then
each point collinear with $U$ and $V$ is also collinear with $W$.

Since $U$ and $V$ are two non-collinear points collinear with $x_p$, we may write, by Lemma~\ref{lemma:coll},
$U=C_p$, $V=D_p$, with $C,D\in\cC$, $x\in C\cap D$, and with $C$ and $D$ not touching each other. The latter
condition implies that $C$ and $D$ share an additional point $y\neq x$. Then $y_p$ is collinear with both $C_p$
and $D_p$. We set $W=E_p$, with $E\in\cC$ and $x\in E$. If $Y=y_p$, then $y\in E$. The points collinear with
$C_p$ and $D_p$ are, besides $x_p$ and $y_p$, all points $F_p$, with $F$ a circle touching both $C$ and $D$.
But by Condition~[CH2], the circle $E$ also touches $F$, so $E_p$ is collinear with $F_p$.

If $Y\neq y_p$, then it is one of the $F_p$ above, and the assertion follows anyway. \qed

Note that the previous proof immediately implies the following lemma.

\begin{lemma}
Every span of $\Gamma$ containing an absolute point of $\rho$ contains exactly two absolute points. Also the
dual holds.
\end{lemma}

In view of the two previous lemmas, there only remains to check Condition~$(iv)$ of Theorem~\ref{main2} in
order to prove that $\Gamma$ is a mixed quadrangle. Therefore we have to look at the dual net corresponding to
a regular point $x_p$, $x \in \mathcal{P}$. In view of the previous results, one can easily give the following
description of the dual net $\Gamma^*_{x_p}$. The points are the circles containing $x$ and the blocks are the
points different from $x$, with incidence given by containment. The circles with gnarl $x$ correspond with a
class of parallel points given by the line $x_l=x_p^\rho$ of the quadrangle $\Gamma$. Then the following
observations are immediate.

\begin{lemma}\begin{itemize}\item[$(i)$] With the above notation, $(\cP,\cC,\partial)$ satisfies Condition~{\rm
[TR]} if and only for each point $x\in\cP$, the dual net $\Gamma^*_{x_p}$ satisfies Axiom~{\rm (LD)} with
respect to the parallel class of points given by the line $x_l$ of $\Gamma$.

\item[$(ii)$] With the above notation, $(\cP,\cC,\partial)$ satisfies Condition~{\rm [F]} if and only if for
each point $x\in\cP$, the dual net $\Gamma^*_{x_p}$ satisfies Axiom~{\rm (VY)}.

\end{itemize}
\end{lemma}

Putting together the last four lemmas, Main Results~\ref{MR:TR} and~\ref{MR:F} follow from Theorem~\ref{main1}
and~\ref{main2}, respectively.

If we substitute Condition~[MP1] with Condition~[MP$1'$], then the dual net $\Gamma^*_{x_p}$ is clearly a dual
affine plane, so Axiom~{\rm (VY)}, or the equivalent Condition~{\rm [F]}, is true trivially. Whence Main
Result~\ref{MR:PR} (the other direction of that theorem being contained in \cite{hvm0}).

Address of the Authors.

Ghent University\\
Department of Pure Mathematics and Computer Algebra\\
Krijgslaan 281, S22,\\
B-9000 Gent\\
Belgium.

\texttt{kstruyve@cage.UGent.be, hvm@cage.UGent.be}

\end{document}